\documentclass[12pt]{article}

\usepackage{amsmath,amssymb,amsthm}
\usepackage {graphicx}

\newtheorem{theorem}{Theorem}[section]
\newtheorem{lemma}[theorem]{Lemma}
\newtheorem{remark}[theorem]{Remark}

\newtheorem{problem}[theorem]{Problem}
\def \RR{\mathbb R}
\def \NN{\mathbb N}
\def\({\left(}  \def\){\right)}
\def\[{\left[}  \def\]{\right]}

\def \beq {\begin {equation}}
\def \eeq {\end{equation}}
\def \OL {\overline}

\def \W {\widetilde}
\def\NN {\rm {\bf N}}

\begin {document}
\large{
\begin {center}
\LARGE{
Absolute Minima of Potentials \\ of Certain Regular Spherical Configurations
}

\bigskip

Sergiy Borodachov

\bigskip

{\it Towson University}

\end {center}

\begin {abstract}
We use methods of approximation theory to find the absolute minima on the sphere of the potential of spherical $(2m-3)$-designs with a non-trivial index $2m$ that are contained in a union of $m$ parallel hyperplanes, $m\geq 2$, whose locations satisfy certain additional assumptions. The interaction between points is described by a function of the dot product, which has positive derivatives of orders $2m-2$, $2m-1$, and $2m$. This includes the case of the classical Coulomb, Riesz, and logarithmic potentials as well as a completely monotone potential of the distance squared. We illustrate this result by showing that the absolute minimum of the potential of the set of vertices of the icosahedron on the unit sphere $S^2$ in $\RR^3$ is attained at the vertices of the dual dodecahedron and the one for the set of vertices of the dodecahedron is attained at the vertices of the dual icosahedron. The absolute minimum of the potential of the configuration of $240$ minimal vectors of $E_8$ root lattice normalized to lie on the unit sphere $S^7$ in $\RR^8$ is attained at a set of $2160$ points on $S^7$ which we describe. 
\end {abstract}



\section {Introduction}

Let 
$$
S^d:=\{(x_1,\ldots,x_{d+1})\subset \RR^{d+1} : x_1^2+\ldots+x_{d+1}^2=1\},
$$
denote the unit sphere in the Euclidean space $\mathbb \RR^{d+1}$, $d\in \NN$. We will call any function
$f:[-1,1]\to(-\infty,\infty]$ continuous on $[-1,1)$ such that $f(1)=\lim\limits_{t\to 1^-}f(t)$ a potential function and will specify additional assumptions on $f$ later. We consider the problem of finding the absolute minimum over the sphere $S^d$ of the potential 
$$
p_f({\bf x},\omega_N):=\sum\limits_{i=1}^{N}f({\bf x}\cdot {\bf x}_i)
$$
of point configurations $\omega_N:=\{{\bf x}_1,\ldots,{\bf x}_N\}$ on $S^d$ satisfying certain regularity assumptions and provide explicit examples for certain dimensions $d$. 
\begin {problem}\label {P1}
{\rm
Given a configuration $\omega_N$ on $S^d$, find the absolute minimum
\begin {equation}\label {P_1}
P_f(\omega_N,S^d):=\min\limits_{{\bf x}\in S^d}\sum\limits_{i=1}^{N}f({\bf x}\cdot {\bf x}_i)
\end {equation}
and points ${\bf x}^\ast\in S^d$ which attain the minimum on the right-hand side of \eqref {P_1}.
}
\end {problem}
We can rewrite the quantity in \eqref {P_1} in the following form
$$
P_f(\omega_N,S^d)=\min\limits_{{\bf x}\in S^d}\sum\limits_{i=1}^{N}g\(\left|{\bf x}-{\bf x}_i\right|^2\),
$$
where $g:[0,4]\to(-\infty,\infty]$ is such that $f(t)=g(2-2t)$. Our greatest interest is in the case when the potential function $f$ is absolutely monotone on $[-1,1)$ or, equivalently, when $g$ is completely monotone on $(0,4]$. Recall that a function $h$ is called {\it absolutely monotone} on an interval $I$ if $f^{(k)}(t)\geq 0$, $t\in I$, $k=0,1,2,\ldots$, and {\it completely monotone} on $I$ if $(-1)^k f^{(k)}(t)\geq 0$, $t\in I$, $k=0,1,2,\ldots$. Furthermore, $h$ is called {\it strictly absolutely monotone} (respectively, {\it strictly completely monotone}) on $I$ if the corresponding inequality is strict in the interior of $I$ for all $k$. Important examples of the kernel $K({\bf x},{\bf y})=g\(\left|{\bf x}-{\bf y}\right|^2\)$ with a (strictly) completely monotone potential function $g$ are the Riesz $s$-kernel, where
\begin {equation}\label {Riesz}
g(t)=\begin {cases}
t^{-s/2}, & s>0,\cr
-t^{-s/2}, & -2< s<0,
\end {cases}
\end {equation}
the logarithmic kernel with $g(t)=\frac {1}{2}\ln \frac {1}{t}$, and the Gaussian kernel
with $g(t)=e^{-\sigma t}$, where $\sigma$ is a positive constant (the Riesz potential function for $-2<s<0$ and the logarithmic potential function become completely monotone on $(0,4]$ after adding an appropriate constant).

The study of the problem about the absolute minimum (maximum) of the potential of a spherical configuration started in mathematical literature several decades ago.
Stolarsky in his early work \cite {Sto1975circle,Sto1975} and later Nikolov and Rafailov \cite {NikRaf2011,NikRaf2013} solved Problem \ref {P1} and the analogous maximization problem for potential function \eqref {Riesz} (and its horizontal translations to the left) with $s\neq 0$ and the vertices of a regular polygon, regular simplex, regular cross-polytope, and a cube. In all these results, for all four types of configurations, the situation is determined by the sign of the derivative of order $n+1$, where $n$ is the strengh of the configuration as a spherical design.

Concerning general potentials, paper \cite {Borsimplex} showed that $p_f({\bf x},\omega_{d+2}^\ast)$, where $\omega^\ast_{d+2}$ is the set of vertices of a regular simplex inscribed in $S^d$, is minimized at points of $-\omega_{d+2}^\ast$ and maximized at points of $\omega_{d+2}^\ast$ for any potential function $f$ with a convex derivative $f'$ on $(-1,1)$.
Problem \ref {P1} was solved in \cite {Borstiff} for $f$ having a convex derivative of order $2m-2$ and configurations $\omega_N$ on $S^d$ called there $m$-stiff. They are $(2m-1)$-designs contained in $m$ parallel hyperplanes. This class includes sets of vertices of a regular $2m$-gon inscribed in $S^1$, of a regular cross-polytope and of the cube inscribed in $S^d$ ($m=2$), of the $24$-cell inscribed in $S^3$ ($m=3$), and many others. The minimum of the potential is attained at every point on $S^d$ that forms $m$ distinct values of the dot product with points of $\omega_N$. The proof of the main result from \cite {Borstiff} can also be found in presentation~\cite {Bor2022talk}.

For sharp configurations (they are $(2m-1)$-designs with $m$ distinct distances between distinct points), both absolute maxima and absolute minima were found when a sharp configuration is a $2m$-design (we call it strongly sharp) and $f$ has a convex derivative of order $2m-1$ (cf. \cite {Borstiff,BorMinMax}). Such are sets of vertices of a regular $(2m+1)$-gon inscribed in $S^1$ and of a regular simplex inscribed in $S^d$ ($m=1$) as well as the Schl\"affi configuration on $S^5$ and the McLaughlin configuration on $S^{21}$ ($m=2$). 

The problem of maxinizing the potential was also solved in \cite {BorMinMax} for antipodal sharp configurations and $f$ having a convex derivative of order $2m-1$. This case, in particular, includes sets of vertices of a regular $2m$-gon inscribed in $S^1$, of a regular cross-polytope inscribed in $S^d$, and of the icosahedron inscribed in $S^2$ as well as the configurations of minimal vectors of $E_8$ root lattice and of the Leech lattice (normalized to lie on $S^7$ and $S^{23}$, respectively). These results inspired further work in \cite {BoyDraHarSafSto600cell}, where the absolute maximum of the potential was found for the $600$-cell on $S^3$ for multiply monotone potentials $f$ (using universal bounds for polarization).

Problems of finding critical points of the total potential of a finite configuration of charges was, in particfular, studied in \cite {Bil2015,GioKhi2020}.

An important application of Problem \ref {P1} is the polarization problem. For a given cardinality $N$, it requires to maximize the minimum value over $S^d$ of the potential of an $N$-point configuration located on $S^d$. It is usually expected that a regular configuration will be optimal (for special cardinalities and/or dimesions). Soluton of Problem \ref {P1} for such configurations would provide a referece constant. Similarly, the maximization counterpart of Problem \ref {P1} 
can be used for the problem of minimizing the maximum value of an $N$-point potential on a sphere. A review on polarization problem can be found, for example, in book \cite {BorHarSafbook} with most recent results reviewed, for example, in \cite {BorMinMax}.

Some sharp configurations are not stiff (for example, icosahedron and minimal non-zero vectors of $E_8$ root lattice). Therefore, the problem of minimizing their potential over the sphere cannot be solved by applying directly the  technique used in \cite {Borstiff} or \cite {BorMinMax}. In this paper, we modify this technique by utilizing non-trivial indexes of these configurations and prove a general result that gives absolute minima of their potentials. This result also applies to some configurations that are not sharp or stiff (e.g., dodecahedron).

In the proof of the main theorem, we utilize the linear programming approach (also known as the Delsarte method or polynomial method), see \cite {DelGoeSei1977}. This technique is widely used among others for best-packing, kissing number, and minimal energy problems, as well as for estimating cardinalities of spherical designs (see, \cite {Lev1979,Lev1992,Lev1998,BoyDanKaz2001,Mus2003,CohKum2007,BorHarSafbook} and references therein).

The paper is structured as follows. We start by mentioning necessary facts and definitions in Section \ref {2}. Section \ref {main} contains the main theorem and its consequences for the particular configurations under consideration. In Section~\ref {Th3.1}, we prove the main result contained in Theorem \ref {generalth}, and in Sections \ref {12}, \ref {20}, and~\ref {240}, we prove its consequences, Theorems \ref {icosahedron}, \ref {dodecahedron}, and~\ref {E_8main}.

\section {Preliminaries on spherical designs}\label {2}

The theory of spherical configurations (codes) and designs is naturally related to the family of Gegenbauer polynomials (cf., e.g., \cite [Chapter 22]{AbrSte1965}).
We define the weight function 
$$
w_d(t):=\gamma_d (1-t^2)^{d/2-1}, 
$$
where the constant $\gamma_d$ is chosen so that $w_d(t)$ is a probability density on $[-1,1]$. Recall that {\it Gegenbauer orthogonal polynomials} (corresponding to~$S^d$) are polynomials in the sequence $\{P_n^{(d)}\}_{n=0}^{\infty}$ such that ${\rm deg}\ \! P_n^{(d)}=n$, $n=0,1,2,\ldots$,  and
$$
\langle P_i,P_j\rangle_d:=\int\limits_{-1}^{1}P_i^{(d)}(t)P_j^{(d)}(t)w_d(t)\ \! dt=0,\ \ i\neq j,
$$ 
normalized so that $P_n^{(d)}(1)=1$, $n=0,1,2,\ldots$.

Given a point configuration $\omega_N=\{{\bf x}_1,\ldots,{\bf x}_N\}\subset S^d$, let $\mathcal I_d(\omega_N)$ be the set of all $n\in \NN$ such that
\begin {equation}\label {Pge}
\sum\limits_{i=1}^{N}\sum\limits_{j=1}^{N}P_n^{(d)}({\bf x}_i\cdot {\bf x}_j)=0.
\end {equation}
We call $\mathcal I_d(\omega_N)$ {\it the index set} of configuration $\omega_N$.
A point configuration $\omega_N=\{{\bf x}_1,\ldots,{\bf x}_N\}\subset S^d$ is called a {\it spherical $m$-design}, $m\in \NN$, if
$$
\frac {1}{N}\sum\limits_{i=1}^{N}p({\bf x}_i)=\int\limits_{S^d}p({\bf x})\ \! d\sigma_d({\bf x})
$$
for every polynomial $p$ on $\RR^{d+1}$ of degree $\leq m$, where $\sigma_d$ is the surface area measure on $S^d$ normalized to be a probability measure. The integer $m$ is called {\it the strength} of spherical design $\omega_N$.
The following well-known statement holds (cf., e.g., \cite [Theorem 5.4.2]{BorHarSafbook}).
\begin {theorem}
A configuration $\omega_N\subset S^d$ is a spherical $m$-design if and only if $\{1,2,\ldots,m\}\subset \mathcal I_d(\omega_N)$.
\end {theorem}
For an antipodal configuration $\omega_N$ on $S^d$, the index set $\mathcal I_d(\omega_N)$ also contains all odd positive integers. 
Following \cite {BoyDanKaz2001}, we call a positive integer $n$ a {\it non-trivial index} of a given antipodal point configuration $\omega_N\subset S^d$ if $n\in \mathcal I_d(\omega_N)$, $n$ is even, and $n>m$, where $m$ is the strength of spherical design $\omega_N$.
Table 1 in the next section (see \cite {BoyDanKaz2001}) lists non-trivial indexes of the configurations considered in this paper. This table can be obtained by verifying directly equality \eqref {Pge}.

%
%
%
%

\section {Main results}\label {main}

For a point ${\bf x}\in S^{d}$, we denote
$$
D({\bf x},\omega_N):=\{{\bf x}\cdot {\bf x}_i : i=\OL {1,N}\}.
$$
The following theorem holds.  

\begin {theorem}\label {generalth}
Let $d,m\in \NN$, $m\geq 2$, and $\omega_N=\{{\bf x}_1,\ldots,{\bf x}_N\}$ be a point configuration on $S^d$ whose index set $\mathcal I_d(\omega_N)$ contains numbers $1,2,\ldots,2m-3,2m-1,2m$. Assume that numbers $-1< t_1<t_2<\ldots<t_m< 1$ are such that 
\begin {equation}\label {sumsq}
\sum_{i=1}^{m}t_i<t_m/2 \ \ \text { and }\ \ \sum\limits_{i=1}^{m}t_i^2-2\(\sum\limits_{i=1}^{m}t_i\)^2<\frac {m(2m-1)}{4m+d-3},
\end {equation}
and that the set $\mathcal D$ of points ${\bf x}^\ast\in S^d$ with $D({\bf x}^\ast,\omega_N)\subset \{t_1,\ldots,t_m\}$ is non-empty.
Let $f:[-1,1]\to (-\infty,\infty]$ be a function continuous on $[-1,1)$ with $\lim\limits_{t\to 1^-}f(t)=f(1)$ and differentiable $2m$ times in $(-1,1)$ with non-negative derivatives $f^{(2m-2)}$, $f^{(2m-1)}$, and $f^{(2m)}$ on $(-1,1)$. Then for every point ${\bf x}^\ast\in \mathcal D$,
\begin {equation}\label {star}
\min\limits_{{\bf x}\in S^d}\sum\limits_{i=1}^{N}f({\bf x}\cdot{\bf x}_i)=\sum\limits_{i=1}^{N}f({\bf x}^\ast\cdot {\bf x}_i).
\end {equation}
If, in addition, $f^{(2m)}>0$ on $(-1,1)$, then the absolute minimum in \eqref {star} is achieved only at points of the set $\mathcal D$. 
\end {theorem}
Statement of Theorem \ref {generalth} was presented in talk \cite {BortalkGMU}. Recall that a point configuration $\omega_N\subset S^d$ is called {\it antipodal} if for every point ${\bf z}\in \omega_N$, the point $-{\bf z}$ is also in $\omega_N$.

\begin {remark}\label {3.2}
{\rm
When $\omega_N$ is antipodal in Theorem \ref {generalth}, conditions \eqref {sumsq} simplify to
\begin {equation}\label {sumsquares}
\sum\limits_{i=1}^{m}t_i^2<\frac {m(2m-1)}{4m+d-3}.
\end {equation}
The quantity on the right-hand side of \eqref {sumsquares} (and in \eqref {sumsq}) is the ratio of the absolute value of the coefficient of $t^{2m-2}$ of Gegenbauer polynomial~$P^{(d)}_{2m}$ to its leading coefficient.
}
\end {remark}

We illustrate Theorem \ref {generalth} by applying it to regular icosahedron and dodecahedron on $S^2$ and to $E_8$ lattice on $S^7$ (see the definition before Theorem~\ref {E_8main}). Table 1 lists deisgn strength and non-trivial indexes for these configurations.

Recall that {\it the regular icosahedron} (denoted by $\OL\omega_{12}$) consists of $12$ points on $S^2$ obtained by cyclic permutations of the coordinates of each of the 4 vectors $(0,\pm1,\pm \varphi)$, where $\varphi=\frac {1+\sqrt{5}}{2}$ is the golden ratio, and normalizing these vectors to lie on $S^2$. 
{\it The regular dodecahedron} (denoted by $\OL\omega_{20}$) is the configuration of $20$ points on $S^2$, where $8$ points have each coordinate $1$ or $-1$ and the remaining $12$ points are obtained by cyclic permutations of the coordinates of the $4$ vectors $\(0,\pm \varphi,\pm 1/\varphi\)$ and normalizing all $20$ points to lie on $S^2$. We will call $\OL\omega_{20}$ {\it dual} for the configuration $\OL\omega_{12}$ and $\OL\omega_{12}$ {\it dual} for the configuration $\OL\omega_{20}$. 

Theorem \ref {generalth} implies the following result for the icosahedron.

\begin {theorem}\label {icosahedron}
Let $f:[-1,1]\to (-\infty,\infty]$ be a function continuous on $[-1,1)$ with $\lim\limits_{t\to 1^-}f(t)=f(1)$ and differentiable $8$ times in $(-1,1)$ with non-negative derivatives $f^{(6)}$, $f^{(7)}$, and $f^{(8)}$ on $(-1,1)$. Then the absolute minimum of the potential 
$$
p_f({\bf x};\OL\omega_{12})=\sum\limits_{{\bf y}\in \OL\omega_{12}}f({\bf x}\cdot {\bf y})
$$
of the vertices of the regular icosahedron $\OL\omega_{12}$ over ${\bf x}\in S^2$ is achieved at every point of the dual regular dodecahedron $\OL \omega_{20}$. If, in addition, $f^{(8)}>0$ on $(-1,1)$, then the absolute minimum of the potential $p_f(\cdot;\OL\omega_{12})$ over $S^2$ is attained only at points of $\OL\omega_{20}$.
\end {theorem}

Theorem \ref {generalth} implies the following result for the dodecahedron.

\begin {theorem}\label {dodecahedron}
Let $f:[-1,1]\to (-\infty,\infty]$ be a function continuous on $[-1,1)$ with $\lim\limits_{t\to 1^-}f(t)=f(1)$ and differentiable $8$ times in $(-1,1)$ with non-negative derivatives $f^{(6)}$, $f^{(7)}$, and $f^{(8)}$ on $(-1,1)$. Then the absolute minimum of the potential 
$$
p_f({\bf x};\OL\omega_{20})=\sum\limits_{{\bf y}\in \OL\omega_{20}}f({\bf x}\cdot {\bf y})
$$
of the vertices of the regular dodecahedron $\OL\omega_{20}$ over ${\bf x}\in S^2$ is achieved at every point of the regular icosahedron $\OL \omega_{12}$. If, in addition, $f^{(8)}>0$ on $(-1,1)$, then the absolute minimum of the potential $p_f(\cdot;\OL\omega_{20})$ over $S^2$ is attained only at points of $\OL\omega_{12}$.
\end {theorem}

Recall that {\it the $E_8$ lattice} is the lattice in $\RR^8$ consisting of those vectors in $\mathbb Z^8 \cup (\mathbb Z+{1}/{2})^8$ with coordinates that sum to an even integer.
We denote by $\OL \omega_{240}$ the set of minimal length non-zero vectors of the $E_8$ lattice normalized to lie on the unit sphere $S^7$. The configuration $\OL \omega_{240}$ consists of $4\(8 \atop 2\)=112$ vectors with $6$ zero coordinates and two coordinates with $\pm 1/\sqrt {2}$ and $2^7=128$ vectors with all eight coordinates $\pm \frac {1}{2\sqrt {2}}$ and even number of ``$-$" signs. For brevity, we will also call $\OL \omega_{240}$ the $E_8$ lattice.

Denote by $\OL \omega_{2160}$ the set of $N=2160$ vectors on $S^7$ that includes $16\(8 \atop 4\)=1120$ vectors with $4$ zero coordinates and $4$ coordinates with $\pm 1/2$, $16$ vectors with $7$ zero coordinates and one coordinate with $\pm 1$, and $8\(16+\(8\atop 3\)+\(8 \atop 5\)\)=1024$ vectors with $7$ coordinates with $\pm 1/4$, one coordinate with $\pm 3/4$, and an odd number of negative coordinates. 

The following result follows from Theorem~\ref {generalth} for the $E_8$ lattice.

\begin {theorem}\label {E_8main}
Let $f:[-1,1]\to (-\infty,\infty]$ be a function continuous on $[-1,1)$ with $\lim\limits_{t\to 1^-}f(t)=f(1)$ and differentiable $10$ times in $(-1,1)$ with non-negative derivatives $f^{(8)}$, $f^{(9)}$, and $f^{(10)}$ on $(-1,1)$. Then the absolute minimum of the potential
$$
p_f({\bf x};\OL\omega_{240})=\sum\limits_{{\bf y}\in \OL\omega_{240}}f({\bf x}\cdot {\bf y})
$$
over ${\bf x}\in S^7$ is achieved at every point of the set $\OL \omega_{2160}$. If, in addition, $f^{(10)}>0$ on $(-1,1)$, then the absolute minimum of the potential $p_f(\cdot;\OL\omega_{240})$ over $S^7$ is attained only at points of $\OL\omega_{2160}$.
\end {theorem}

\begin {table}\label {table1}
\centering
\begin{tabular}{ |c|c|c| } 
 \hline
Name or notation & Strength as a spherical design & Non-trivial indexes \\ 
\hline
 Icosahedron $\OL \omega_{12}$& $5 $& $8$, $14$ \\ 
 Dodecahedron $\OL\omega_{20}$& $5$ & $8$, $14$ \\ 

$E_8$ lattice $\OL\omega_{240}$& $7$ & $10$ \\
\hline

\end{tabular}
\caption{Non-trivial indexes of spherical configurations considered in the paper}
\end {table}

%
%
%

\section {Proof of Theorem \ref {generalth}}\label {Th3.1}

Gegenbauer polynomials are known to satisfy the following recurrence relation (cf., e.g., \cite [Theorem 5.3.1]{BorHarSafbook})
\begin {equation}\label {RR}
P_{n+1}^{(d)}(t)=\frac {2n+d-1}{n+d-1}tP_n^{(d)}(t)-\frac {n}{n+d-1}P_{n-1}^{(d)}(t),\ \ n\geq 1,
\end {equation}
where $P_0^{(d)}(t)=1$ and $P_1^{(d)}(t)=t$. Relation \eqref {RR} implies that for $n$ even, $P_n^{(d)}(t)$ has only even powers of $t$ and for $n$ odd, only odd powers of $t$.
Therefore, we use the notation 
$$
P_n^{(d)}(t)=\alpha_nt^n-\gamma_nt^{n-2}+...,\ \ n\geq 2.
$$
We also let $\alpha_0=\alpha_1=1$. 
The following auxiliary statement holds.
\begin {lemma}\label {n(n-1)}
For every $d\geq 1$ and $n\geq 2$, we have $\alpha_n,\gamma_n>0$ and
\begin {equation}\label {r_n}
r_n:=\frac {\gamma_n}{\alpha_n}=\frac {n(n-1)}{2(2n+d-3)}.
\end {equation}
\end {lemma}
\begin {proof}
Since from \eqref {RR} we have $\alpha_{n+1}=\frac {2n+d-1}{n+d-1}\alpha_n$, $n\geq 0$, the positiveness of the leading coefficients $\alpha_n$ follows by induction. Relation \eqref {RR} also implies that $\gamma_2=\frac {1}{d}$ and that
\begin {equation}\label {rel}
\gamma_{n+1}=\frac {2n+d-1}{n+d-1}\gamma_n+\frac {n}{n+d-1}\alpha_{n-1},\ \ n\geq 2.
\end {equation}
Using \eqref {rel}, one can now show by induction that $\gamma_n>0$, $n\geq 2$. Dividing both sides of \eqref {rel} by $\alpha_{n+1}$, we obtain
$$
r_{n+1}=\frac {2n+d-1}{n+d-1}\cdot \frac {r_n\alpha_n}{\alpha_{n+1}}+\frac {n}{n+d-1}\cdot \frac {\alpha_{n-1}}{\alpha_{n+1}}=r_n+\frac {n(n+d-2)}{(2n+d-3)(2n+d-1)}
$$
for $n\geq 2$ with $r_2=\gamma_2/\alpha_2=\frac {1}{d+1}$. Then \eqref {r_n} follows by induction.
\end {proof}


Denote by $\mathbb P_n$, $n\geq 0$, the set of all polynomials of degree at most $n$.
\begin {lemma}\label {Q}
Let $m\geq 2$ and $-1\leq t_1<\ldots<t_m\leq 1$ be arbitrary real numbers satisfying \eqref {sumsq}.
Then
$$
\langle P_{2m-2}^{(d)},Q_i\rangle_d>0,\ \ i=0,1,2,
$$
where $Q_i(t):=(t-t_m)^{i}\prod\limits_{j=1}^{m-1}(t-t_j)^2$.
\end {lemma}
\begin {proof}
Let $u(t):=t^{2m-2}$. Since $P_{2m-2}^{(d)}\bot \mathbb P_{2m-3}$, we have
\begin {equation}\label {a_m}
\langle P_{2m-2}^{(d)},u\rangle_d=\frac {\langle P_{2m-2}^{(d)},P_{2m-2}^{(d)}\rangle_d}{\alpha_{2m-2}}>0.
\end {equation}
Hence,
$$
\langle P_{2m-2}^{(d)},Q_0\rangle_d=\langle P_{2m-2}^{(d)},u\rangle_d>0.
$$
Observe now that $
Q_1(t)=t^{2m-1}+ct^{2m-2}+\ldots ,
$
where $c=-\(t_m+2\sum_{i=1}^{m-1}t_i\)$. In view of \eqref {sumsq}, we have $c>0$.
Since $P_{2m-2}^{(d)}\bot \mathbb P_{2m-3}$ and $P_{2m-2}^{(d)}\bot P$, where $P(t)=t^{2m-1}$, using \eqref {a_m}, we have 
$$
\langle P_{2m-2}^{(d)},Q_1\rangle_d=c\langle P_{2m-2}^{(d)},u\rangle_d>0.
$$
Denote by $\tau_1,\ldots,\tau_{2m}$ the elements of the multiset $\{t_1,t_1,t_2,t_2,\ldots,t_m,t_m\}$. Then 
$$
Q_2(t)=t^{2m}-at^{2m-1}+bt^{2m-2}+\ldots ,
$$
where $a=\sum\limits_{i=1}^{2m}\tau_i$ and
\begin {equation*}
\begin {split}
b=\sum\limits_{1\leq i<j\leq 2m}\tau_i\tau_j&=\frac {1}{2}\sum\limits_{1\leq i\neq j\leq 2m}\tau_i\tau_j=\frac {1}{2}\sum\limits_{i=1}^{2m}\sum\limits_{j=1}^{2m}\tau_i\tau_j-\frac {1}{2}\sum\limits_{i=1}^{2m}\tau_i^2\\
&=\frac {1}{2}\(\sum\limits_{i=1}^{2m}\tau_i\)^2-\sum\limits_{i=1}^{m}t_i^2=2\(\sum\limits_{i=1}^{m}t_i\)^2-\sum\limits_{i=1}^{m}t_i^2.
\end {split}
\end {equation*}
The assumption in \eqref {sumsq} and Lemma \ref {n(n-1)} impliy that 
$$
\epsilon:=\frac{\gamma_{2m}}{\alpha_{2m}}+b=\frac {m(2m-1)}{4m+d-3}+b>0.
$$ 
Using the monic Gegenbauer polynomial $\W P_{2m}^{(d)}(t)=t^{2m}-\frac {\gamma_{2m}}{\alpha_{2m}}t^{2m-2}+\ldots$, we will have
$$
Q_2(t)=t^{2m}-at^{2m-1}-\frac {\gamma_{2m}}{\alpha_{2m}}t^{2m-2}+\epsilon t^{2m-2}+...=\W P_{2m}^{(d)}(t)-at^{2m-1}+\epsilon u(t)+P(t),
$$
where $P\in \mathbb P_{2m-3}$ is some polynomial. Taking into account \eqref {a_m}, we then have
$$
\langle P_{2m-2}^{(d)},Q_2\rangle_d=\epsilon\langle P_{2m-2}^{(d)},u\rangle_d>0,
$$
which completes the proof of the lemma.
\end {proof}

We next state a known result from numerical analysis (cf., e.g., \cite {IsaKel1965}) on Newton-type interpolation. For completeness, we present its proof.
We define $\Pi_0(t):=1$ and let 
$$
\Pi_i(t):=(t-t_1)^2\cdot\ldots\cdot (t-t_{i})^2, \ \ \ i=1,\ldots,m.
$$

\begin {lemma}\label {interp}
Let $f$ be a function continuous on an interval $I$ and differentiable in its interior. Let $t_1,\ldots,t_m\in {\rm int}\ \! I$, $m\geq 2$, be arbitrary pairwise distinct nodes. Then
\begin {enumerate}
\item [(i)]there exists a vector ${\bf a}=(a_1,b_1,\ldots,a_m,b_m)\in \RR^{2m}$ such that the polynomial
\begin {equation}\label {a}
p(t;{\bf a}):=\sum\limits_{i=1}^{m}(a_i+b_i(t-t_i))\Pi_{i-1}(t)
\end {equation}
satisfies $p(t_i;{\bf a})=f(t_i)$ and $p'(t_i;{\bf a})=f'(t_i)$ for $i=1,\ldots,m$;

\item [(ii)] the vector ${\bf a}$ in (i) is unique and $a_1=f(t_1)$ and $b_1=f'(t_1)$.

Assume, in addition, that $f$ is $2m$ times differentiable in the interior of $I$. Then 

\item [(iii)] for every $i=2,\ldots,m$, there are points $c_i,d_i\in (t_1,t_i)$ such that 
$$
a_i=\frac{f^{(2i-2)}(c_i)}{(2i-2)!}\text{  and  } b_i=\frac {f^{(2i-1)}(d_i)}{(2i-1)!};
$$

\item [(iv)] if, in addition, $f^{(2m)}$ is non-negative in the interior of $I$, we have
\begin {equation}\label {leq}
p(t;{\bf a})\leq f(t),\ \ t\in I,
\end {equation}
and if $f^{(2m)}>0$ in the interior of $I$, then \eqref {leq} is strict for $t\in I\setminus\{ t_1,\ldots,t_m\}$.
\end {enumerate}
\end {lemma}
\begin {proof}
The system of polynomials $\mathcal Q:=\{p_0,p_1,\ldots,p_{2m-1}\}$, where $p_{2i}(t)=\Pi_i(t)$, $i=0,1,\ldots,m-1$, and $p_{2i-1}(t)=(t-t_i)\Pi_{i-1}(t)$, $i=1,\ldots,m$, satisfies ${\rm deg}\ \! p_j=j$, $j=0,1,\ldots,2m-1$. Hence, it forms a basis for $\mathbb P_{2m-1}$. Let $h$ be the Hermite interpolating polynomial for $f$ at $t_1,\ldots,t_m$. Since $h\in \mathbb P_{2m-1}$, there is a vector ${\bf a}:=(a_1,b_1,\ldots,a_m,b_m)\in \mathbb R^{2m}$ such that
$$
h(t)=\sum\limits_{i=1}^{m}a_ip_{2(i-1)}(t)+\sum\limits_{i=1}^{m}b_ip_{2i-1}(t)=\sum\limits_{i=1}^{m}(a_i+b_i(t-t_i))\Pi_{i-1}(t),
$$
which proves (i). The uniqueness of vector ${\bf a}$ follows from the uniqueness of the Hermite interpolating polynomial $h$ and of the coefficients relative to basis $\mathcal Q$. Since $f(t_1)=h(t_1)=a_1$ and $f'(t_1)=h'(t_1)=b_1$, we have (ii). We let $p(t;{\bf a}):=h(t)$.

Observe that for every $k=1,\ldots,m$, the polynomial 
$$
q_k(t):=\sum_{i=1}^{k}(a_i+b_i(t-t_i))\Pi_{i-1}(t)
$$
with $a_i$'s and $b_i$'s as in (i) interpolates $f$ and $f'$ at $t_1,\ldots,t_k$. This occurs because the remainder $v(t):=\sum\limits_{i=k+1}^{m}(a_i+b_i(t-t_i))\Pi_{i-1}(t)$ vanishes together with its derivative at $t_1,\ldots,t_k$.

Choose any index $2\leq k\leq m$ and assume that $f^{(2m)}$ exists in ${\rm int}\ \! I$. By the Rolle's theorem, the difference $f'-q'_k$ has at least $2k-1$ zeros on $[t_1,t_k]$. Applying the Rolle's theorem, $2k-2$ more times, one can find a point $d_k\in (t_1,t_k)$ such that $(f-q_k)^{(2k-1)}(d_k)=0$. Since $q_k^{(2k-1)}(t)=(2k-1)!b_k$, we have $b_k=f^{(2k-1)}(d_k)/(2k-1)!$. The polynomial 
$$
w(t):=q_{k-1}(t)+a_k\Pi_{k-1}(t)=q_k(t)-b_k(t-t_k)\Pi_{k-1}(t)
$$
interpolates $f$ at $t_1,\ldots,t_k$ and $f'$ at $t_1,\ldots,t_{k-1}$.   By the Rolle's theorem, the derivative $(f-w)'$ has at least $2k-2$ zeros in $[t_1,t_k)$, and the derivative $(f-w)^{(2k-2)}$ has at least one zero in $(t_1,t_k)$. That is, there exists $c_k\in (t_1,t_k)$ such that $w^{(2k-2)}(c_k)=f^{(2k-2)}(c_k)$. Since $w^{(2k-2)}(t)=(2k-2)!a_k$, we have $a_k=f^{(2k-2)}(c_k)/(2k-2)!$, which proves (iii).

Let now $\tau\in I\setminus \{t_1,\ldots,t_m\}$ be arbitrary. Define $q_{m+1}(t):=p(t;{\bf a})+a_{m+1}\Pi_{m}(t)$, where $a_{m+1}$ is an arbitrary constant. Then $q_{m+1}$ interpolates $f$ and $f'$ at $t_1,\ldots,t_m$. We choose $a_{m+1}$ such that $q_{m+1}(\tau)=f(\tau)$. By the Rolle's theorem, the difference $f'-q'_{m+1}$ vanishes at $2m$ distinct points of~${\rm int}\ \!I$. Applying the Rolle's theorem $2m-1$ more times, we obtain that the difference $f^{(2m)}-q_{m+1}^{(2m)}$ vanishes at least once (at some point $z=z(\tau)\in {\rm int} \ \! I$). Since $q_{m+1}^{(2m)}(t)\equiv (2m)!a_{m+1}$, we have $a_{m+1}=f^{(2m)}(z)/(2m)!$. Then, since $f^{(2m)}\geq 0$ in ${\rm int}\ \! I$, we obtain that
\begin {equation}\label {tau}
f(\tau)=q_{m+1}(\tau)=p(\tau;{\bf a})+\frac {f^{(2m)}(z)}{(2m)!}\Pi_{m}(\tau)\geq p(\tau;{\bf a})
\end {equation}
for $\tau\neq t_1,\ldots,t_m$ (for $\tau=t_1,\ldots,t_m$, we have equality in \eqref {leq}). If $f^{(2m)}>0$ in ${\rm int}\ \! I$, then the inequality in \eqref {tau} is strict.
\end {proof}

\begin {lemma}\label {perpendicular}
Suppose the nodes $-1< t_1<\ldots< t_m<1$, $m\geq 2$, satisfy conditions~\eqref {sumsq}. Let $f:[-1,1]\to (-\infty,\infty]$ be continuous on $[-1,1)$ with $\lim\limits_{t\to 1^-}f(t)=f(1)$ and differentiable $2m$ times on $(-1,1)$ with non-negative derivatives $f^{(2m-2)}$, $f^{(2m-1)}$, and $f^{(2m)}$ on $(-1,1)$. 

Then for every $d\in \NN$, there is a polynomial $\lambda\in \mathbb P_{2m}$ such that $\lambda\bot P_{2m-2}^{(d)}$, $\lambda(t_i)=f(t_i)$, $i=1,\ldots,m$, and $\lambda(t)\leq f(t)$, $t\in [-1,1]$. If, in addition, $f^{(2m)}>0$ on $(-1,1)$, then $\lambda(t)<f(t)$, $t\in [-1,1]\setminus\{t_1,\ldots,t_m\}$.
\end {lemma}

\begin {proof}
By Lemma \ref {interp}, there exists a unique vector ${\bf a}=(a_1,b_1,\ldots,a_m,b_m)\in \RR^{2m}$ such that the polynomial 
$$
p(t;{\bf a})=\sum\limits_{i=1}^{m}(a_i+b_i(t-t_i))\Pi_{i-1}(t)
$$ 
defined by \eqref {a} interpolates $f$ and $f'$ at the nodes $t_1,\ldots,t_m$ and 
\begin {equation}\label {3a}
p(t;{\bf a})\leq f(t),\ \ \ t\in [-1,1).
\end {equation}
Furthermore, if $f(1)<\infty$, then \eqref {3a} holds at $t=1$ by continuity and, if $f(1)=\infty$, then \eqref {3a} holds at $t=1$ trivially.

Let $\beta_0$ and $\beta_1$ be constants such that $Q_i+\beta_i Q_2\bot P_{2m-2}^{(d)}$, $i=0,1$ (see the notation in the statement of Lemma~\ref {Q}). Lemma \ref {Q} implies that $\langle Q_i,P_{2m-2}^{(d)}\rangle_d>0$, $i=0,1,2$. Hence, $\beta_i<0$, $i=0,1$. Define
\begin {equation*}
\begin {split}
\lambda(t):&=
\sum\limits_{i=1}^{m-1}(a_i+b_i(t-t_i))\Pi_{i-1}(t)+a_m(Q_0(t)+\beta_0Q_2(t))+b_{m}(Q_1(t)+\beta_1Q_2(t))\\
&=\sum\limits_{i=1}^{m}(a_i+b_i(t-t_i))\Pi_{i-1}(t)
+\(a_m\beta_0+b_{m}\beta_1\)\Pi_{m}(t).
\end {split}
\end {equation*}
The polynomial $\lambda\in \mathbb P_{2m}$ interpolates $f$ and $f'$ at $t_1,\ldots,t_m$ and $\lambda\bot P_{2m-2}^{(d)}$. Lemma~\ref {interp} and assumptions of Lemma \ref {perpendicular} imply that $a_m=\frac {f^{(2m-2)}(c_m)}{(2m-2)!}\geq 0$ and $b_m=\frac {f^{(2m-1)}(d_m)}{(2m-1)!}\geq 0$, where $c_m$ and $d_m$ are some points in $(-1,1)$. Since $\Pi_{m}$ is a non-negative polynomial, we obtain that
$
\lambda(t)\leq p(t;{\bf a})\leq f(t),\ t\in [-1,1]$. 

If, in addition, $f^{(2m)}>0$ on $(-1,1)$, then Lemma \ref {interp} implies that $\lambda(t)\leq p(t;{\bf a})< f(t)$, $t\in [-1,1)\setminus\{t_1,\ldots,t_m\}$. If $f(1)<\infty$, then Lemma \ref {interp} can be used with $I=[-1,1]$ to conclude also that $\lambda(1)\leq p(1;{\bf a})< f(1)$. If $f(1)=\infty$, then we conclude this trivially.
\end {proof}

We will also need the following known statement, see \cite {BoyDanKaz2001} (its proof can also be found in \cite [Lemma 5.2.2]{BorHarSafbook}).
\begin {lemma}\label {design}
Let $d,n\in \NN$ and $\omega_N=\{{\bf x}_1,\ldots,{\bf x}_N\}\subset S^d$ be an arbitrary configuration. Then
$$
\sum\limits_{i=1}^{N}\sum\limits_{j=1}^{N}P_n^{(d)}({\bf x}_i\cdot {\bf x}_j)=0\ \ \text {if and only if}\ \ \sum\limits_{i=1}^{N}P_n^{(d)}({\bf x}\cdot {\bf x}_i)=0\ \ \text {for all}\ {\bf x}\in S^d.
$$
\end {lemma}
We are now ready to prove Theorem \ref {generalth}.

\begin {proof}[Proof of Theorem \ref {generalth}]
By Lemma \ref {perpendicular}, there is a polynomial $\lambda\in \mathbb P_{2m}$ such that $\lambda\bot P_{2m-2}^{(d)}$, $\lambda(t_i)=f(t_i)$, $i=1,\ldots,m$, and $\lambda(t)\leq f(t)$, $t\in [-1,1]$. 

Let ${\bf x}^\ast$ be any point in the set $\mathcal D$. Then $\lambda({\bf x}^\ast\cdot {\bf x}_i)=f({\bf x}^\ast\cdot {\bf x}_i)$, $i=1,\ldots,N$. Let $\beta_i$, $i=0,1,\ldots,2m-3,2m-1,2m$, be such that $\lambda(t)=\sum\limits_{n=0\atop n\neq 2m-2}^{2m}\beta_nP_n^{(d)}(t)$. Since the index set of $\omega_N$ contains numbers $1,\ldots,2m-3,2m-1,2m$, by Lemma~\ref 
{design}, for every ${\bf x}\in S^d$, we have 
$$
\sum_{i=1}^{N}P_n^{(d)}({\bf x}\cdot {\bf x}_i)=0,\ \ n=1,\ldots,2m-3,2m-1,2m,
$$ 
and
$$
\sum\limits_{i=1}^{N}\lambda({\bf x}\cdot{\bf x}_i)=\sum\limits_{n=0\atop n\neq 2m-2}^{2m}\beta_n\sum\limits_{i=1}^{N}P_n^{(d)}({\bf x}\cdot {\bf x}_i)=\beta_0\sum\limits_{i=1}^{N}P_0^{(d)}({\bf x}\cdot{\bf x}_i)=N\beta_0.
$$
Then for every ${\bf x}\in S^d$, we have
\begin {equation}\label {q1b}
\sum\limits_{i=1}^{N}f({\bf x}\cdot {\bf x}_i)\geq \sum\limits_{i=1}^{N}\lambda({\bf x}\cdot {\bf x}_i)=N\beta_0=\sum\limits_{i=1}^{N}\lambda({\bf x}^\ast\cdot {\bf x}_i)=\sum\limits_{i=1}^{N}f({\bf x}^\ast\cdot {\bf x}_i),
\end {equation}
which implies \eqref {star}.

Assume now that $f^{(2m)}>0$ on $(-1,1)$ and let ${\bf x}\in S^d\setminus \mathcal D$ be any point. Then $D({\bf x},\omega_N)\not\subset \{t_1,\ldots,t_m\}$. Then for some $1\leq j\leq N$, we have ${\bf x}\cdot {\bf x}_j\neq t_1,\ldots,t_m$. Lemma \ref {perpendicular} now implies that $\lambda({\bf x}\cdot{\bf x}_j)<f({\bf x}\cdot{\bf x}_j)$, which makes the inequality in \eqref {q1b} strict. Hence, the minimum in \eqref {star} is not attained at the point ${\bf x}$.
\end {proof}

\section {Proof of Theorem \ref {icosahedron}}\label {12}

We will first establish two auxiliary statements.
\begin {lemma}\label {psi}
Suppose ${\bf x}=(x_1,x_2,x_3)\in \RR^3$ is a vector with non-negative coordinates and the maximal coordinate equal to $1$. Then ${\bf x}$ forms at least four distinct dot prducts with $12$ vectors of the set 
$$
A:=\{(0,\pm 1,\pm\varphi),(\pm 1,\pm\varphi,0),(\pm \varphi,0,\pm 1)\}. 
$$
Vector ${\bf x}$ forms exactly four distinct dot products with vectors from $A$ if and only if one of the following holds
\begin {enumerate}
\item [(i)] ${\bf x}\in A_1:=\{(0,1/\varphi,1),(1/\varphi,1,0),(1,0,1/\varphi)\}$ with dot products $\pm 1,\pm \sqrt {5}$

\item [(ii)] or ${\bf x}\in A_2:=\{(0,1,1/\varphi^2),
(1/\varphi^2,0,1),(1,1/\varphi^2,0)\}$ with dot products $\pm 1\pm 1/\varphi$

\item [(iii)] or ${\bf x}=(1,1,1)$ with dot products $\pm 1\pm \varphi$.
\end {enumerate}
\end {lemma}
\begin {proof}
Assume that ${\bf x}$ forms at most four dot products with vectors from $A$.
If $x_1=1$, we let $b=x_2$ and $c=x_3$, if $x_1<1$ and $x_2=1$, we let $b=x_3$ and $c=x_1$. If $x_1,x_2<1$ then $x_3=1$ and we let $b=x_1$ and $c=x_2$. Then we have $0\leq b,c\leq 1$ and ${\bf x}$ forms all of the following dot products with vectors from $A$: $\pm \varphi\pm c$, $\pm b\pm c\varphi$, and $\pm 1\pm b\varphi$. 

Assume first that $0<c<1$. Since $\varphi>c>0$, the set $G=\{\pm \varphi\pm c\}$ contains four distinct dot products and, hence, dot products $1\pm b\varphi$ must be in $G$.
Observe that $1+b\varphi\geq 1-b\varphi>-\varphi+c$. Assume to the contrary that $b>0$. Then $1-b\varphi=\varphi-c$ and $1+b\varphi=\varphi+c$. Adding these two equations we obtain a contradiction $\varphi=1$. Therefore, $b=0$. Then $1\pm b\varphi=1<\varphi+c$. Therefore, $1=\varphi-c$, which implies that 
$c=\varphi-1=1/\varphi$.  Then ${\bf x}\in A_1$.

Asume next that $c=0$. Then ${\bf x}$ forms the following dot products with vectors from $A$: $\pm \varphi$, $\pm b$, and $\pm 1\pm b\varphi$. If $b=0$, then there are five distinct dot products: $\pm \varphi$, $\pm 1$, and $0$ contradicting the choice of ${\bf x}$. Then $0<b\leq 1$ and the four dot products $\pm \varphi$ and $\pm b$ are all distinct. Since no more dot products can be formed, we have $1+b\varphi=\varphi$ because $1+b\varphi$ is greater than the other three dot products. Consequently, $b=1-1/\varphi=1/\varphi^2$ and ${\bf x}\in A_2$. 

Finally, assume that $c=1$. Then some of the dot products are $-\varphi-1<-\varphi+1<\varphi-1<\varphi+b\leq \varphi+1$. We must have $b=1$. Consequently, ${\bf x}=(1,1,1)$.

Thus, if ${\bf x}\notin A_1\cup A_2\cup \{(1,1,1)\}$ then ${\bf x}$ forms more than four distinct dot products with vectors from $A$. If ${\bf x}\in A_1\cup A_2\cup \{(1,1,1)\}$ then ${\bf x}$ forms exactly four distinct dot products with vectors from $A$. They are listed in (i), (ii), and~(iii).
\end {proof}

\begin {lemma}\label {icosahdr}
For every ${\bf u}\in S^2$, the set $D({\bf u},\OL\omega_{12})$ has at least four distinct elements. The only vectors ${\bf u}\in S^2$ with $D({\bf u},\OL\omega_{12})$ having exactly four distinct dot products are those in $\OL\omega_{20}$ (the dot products are $\frac {\pm \varphi\pm 1}{\sqrt {3\varphi+6}}$) and in $\OL\omega_{12}$ itself (the dot products are $\pm 1$ and $\pm \frac {1}{\sqrt {5}}$).
\end {lemma}
\begin {proof}
For any vector ${\bf u}=(u_1,u_2,u_3)\in S^2$, let ${\bf x}=\frac {1}{a}{\bf y}$, where ${\bf y}=(\left|u_1\right|,\left|u_2\right|,\left|u_3\right|)$ and $a=\max\{\left|u_i\right| : i=\OL {1,3}\}$. We have $a>0$, since ${\bf u}\neq {\bf 0}$. The set $D({\bf u},\OL\omega_{12})$, where we have $\OL\omega_{12}=\frac {1}{\sqrt {\varphi+2}} A$, equals the set of dot products that ${\bf x}$ forms with vectors from $\frac {a}{\sqrt {\varphi+2}}A$. By Lemma \ref {psi}, vector ${\bf u}$ must form at least four distinct dot products with vectors from $\OL\omega_{12}$. Vector ${\bf u}$ forms exactly four distinct dot products if and only if ${\bf x}$ is in $A_1$ or in $A_2$ or equals  $(1,1,1)$. We have ${\bf u}\in \OL\omega_{12}$ if and only if ${\bf x}\in A_1$ and ${\bf u}\in \OL\omega_{20}$ if and only if ${\bf x}\in A_2$ or ${\bf x}=(1,1,1)$. 

If ${\bf u}\in \OL\omega_{12}$, then $a=\frac {\varphi}{\sqrt {\varphi+2}}$ and $\frac {a}{\sqrt {\varphi+2}}=\frac {\varphi}{\varphi+2}=\frac {1}{\sqrt {5}}$ and ${\bf x}\in A_1$. In view of Lemma \ref {psi} (i),
$$
D({\bf u},\OL\omega_{12})=D\({\bf x},\frac {1}{\sqrt{5}}A\)=\left\{\pm \frac {1}{\sqrt {5}},\pm 1\right\}.
$$
In the case ${\bf u}\in \OL\omega_{20}$, we have ${\bf x}\in A_2$ or ${\bf x}=(1,1,1)$. If ${\bf x}\in A_2$, then $a=\frac {\varphi}{\sqrt {3}}$ and $\frac {a}{\sqrt {\varphi+2}}=\frac {\varphi}{\sqrt {3\varphi+6}}$. In view of Lemma \ref {psi} (ii), $D({\bf u},\OL\omega_{12})=D\({\bf x},\frac {\varphi}{\sqrt {3\varphi+6}}A\)=\left\{\frac {\pm \varphi\pm 1}{\sqrt {3\varphi+6}}\right\}$. Finally, if ${\bf x}=(1,1,1)$, then $a=1/\sqrt {3}$ and $\frac {a}{\sqrt {\varphi+2}}=\frac {1}{\sqrt {3\varphi+6}}$. Lemma \ref {psi} (iii) then yields $D({\bf u},\OL\omega_{12})=\left\{\frac {\pm \varphi\pm 1}{\sqrt {3\varphi+6}}\right\}$.
\end {proof}

\begin {proof}[Proof of Theorem \ref {icosahedron}]The icosahedron $\OL\omega_{12}=\{{\bf x}_1,\ldots,{\bf x}_{12}\}$ is an antipodal configuration on $S^2$ whose index set $\mathcal I_2(\OL\omega_{12})$ contains the set $\{1,2,3,4,5,7,8\}$, see Table 1. By Lemma \ref {icosahdr}, any point ${\bf x}^\ast\in \OL\omega_{20}$ has $m=4$ distinct dot products with vectors from $\OL\omega_{12}$, which are $t_1=-\frac {\varphi+1}{\sqrt {3\varphi+6}}$, $t_2=\frac {-\varphi+1}{\sqrt {3\varphi+6}}$, $t_3=\frac {\varphi-1}{\sqrt {3\varphi+6}}$, $t_4=\frac {\varphi+ 1}{\sqrt {3\varphi+6}}$. Observe that
$$
\sum\limits_{i=1}^{4}t_i^2=\frac {2(\varphi+1)^2}{3\varphi+6}+\frac {2(\varphi-1)^2}{3\varphi+6}=\frac {4}{3}<\frac {28}{15}=\frac {m(2m-1)}{4m+d-3},
$$
where $m=4$. In view of Remark \ref {3.2}, conditions \eqref {sumsq} hold. Then by Theorem~\ref {generalth},
$$
\min\limits_{{\bf x}\in S^2}\sum\limits_{i=1}^{12}f({\bf x}\cdot {\bf x}_i)=\sum\limits_{i=1}^{12}f({\bf x}^\ast\cdot {\bf x}_i)
$$
for any point ${\bf x}^\ast\in \OL\omega_{20}$.

Lemma \ref {icosahdr} implies that for every ${\bf x}\in S^2\setminus \OL\omega_{20}$, the set $D({\bf x},\OL\omega_{12})$ is not contained in $\{\frac {\pm\varphi\pm 1}{\sqrt {3\varphi+6}}\}$. Then by Theorem \ref {generalth}, if $f^{(8)}>0$ on $(-1,1)$, point ${\bf x}$ is not a point of absolute minimum of $p_f(\cdot;\OL\omega_{12})$ over $S^2$.
\end {proof}

\section {Proof of Theorem \ref {dodecahedron}}\label {20}

We start by establishing two auxiliary statements similar to Lemmas \ref {psi} and~\ref {icosahdr}.
\begin {lemma}\label {dodecahedrn}
Suppose ${\bf x}=(x_1,x_2,x_3)\in \RR^3$ is a vector with non-negative coordinates and the maximal coordinate equal to $1$. Then ${\bf x}$ forms at least four distinct dot prducts with $20$ vectors of the set 
$$
B:=\{(0,\pm\varphi,\pm 1/\varphi),(\pm\varphi,\pm 1/\varphi,0),(\pm 1/\varphi,0,\pm \varphi),(\pm 1,\pm 1,\pm 1)\}. 
$$
Vector ${\bf x}$ forms exactly four distinct dot products with vectors from $B$ if and only if 
${\bf x}\in A_1=\{(0,1/\varphi,1),(1/\varphi,1,0),(1,0,1/\varphi)\}$ with the set of dot products being $\{\pm 1\pm 1/\varphi\}=\{\pm \varphi,\pm 1/\varphi^2\}$.
\end {lemma}
\begin {proof}
Assume that ${\bf x}$ forms at most four distinct dot products with vectors from $B$.
If $x_1=1$, we denote $b=x_2$ and $c=x_3$. If $x_1<1$ and $x_2=1$, we let $b=x_3$ and $c=x_1$. Finally, if $x_1<1$ and $x_2<1$, then $x_3=1$ and we let $b=x_1$ and $c=x_2$. Then we have $0\leq b,c\leq 1$ and ${\bf x}$ forms the following dot products with vectors from $B$: $\pm b\varphi\pm c/\varphi$, $\pm \varphi\pm b/\varphi$, $\pm 1/\varphi\pm c\varphi$, and $\pm 1\pm b\pm c$.

Assume to the contrary that $b>0$. Since we also have $b\leq 1$, the set $H:=\{\pm \varphi\pm b/\varphi\}$ contains four distinct dot products, and, hence, the dot product $b\varphi-c/\varphi$ must be in~$H$. Observe that 
$$
b\varphi-c/\varphi<\varphi+b/\varphi\ \ \text {and} \ \ b\varphi-c/\varphi>-\varphi+b/\varphi.
$$
Then we must have $b\varphi-c/\varphi=\varphi-b/\varphi$. Similarly, the dot product $1/\varphi-c\varphi$ must also be in $H$. Since
$$
\varphi-b/\varphi=\varphi-1/\varphi+{(1-b)}/{\varphi}=1+{(1-b)}/{\varphi}>{1}/{\varphi}-c\varphi,
$$
and $-\varphi-b/\varphi<1/\varphi-c\varphi$, we must have $-\varphi+b/\varphi=1/\varphi-c\varphi$. Solving the system of these two equations we obtain that $b=c=1$. Then ${\bf x}=(1,1,1)$. However, vector $(1,1,1)$ forms six distinct dot products with vectors from $B$: $\pm 1=\pm (\varphi-1/\varphi)$, $\pm \sqrt {5}=\pm (\varphi+1/\varphi)$, and $\pm 3$ contradicting the original choice of ${\bf x}$. 

Thus, we must have $b=0$. Then ${\bf x}$ forms the following dot products with vectors from $B$: $\pm c/\varphi$, $\pm \varphi$, $\pm 1/\varphi\pm c\varphi$, and $\pm 1\pm c$. If $c=0$ or $c=1$, then there are more than four distinct dot products. If $0<c<1$, then the set $\{\pm 1\pm c\}$ contains four distinct dot products and $\varphi$ must be in this set. Consequently, $\varphi=1+c$, since $\varphi$ is greater than three other dot products in this set. Therefore, $c=\varphi-1=1/\varphi$; that is, ${\bf x}$ must be in $A_1$. 

Thus, if ${\bf x}\notin A_1$, then ${\bf x}$ forms more that four distinct dot products with vectors from $B$. If ${\bf x}\in A_1$, then ${\bf x}$ forms four distinct dot products with vectors from $B$: $\pm (1+1/\varphi)=\pm \varphi$ and $\pm (1-1/\varphi)=\pm 1/\varphi^2$.
\end {proof}

\begin {lemma}\label {ddcd}
For every ${\bf u}\in S^2$, the set $D({\bf u},\OL\omega_{20})$ has at least four distinct elements. The only vectors ${\bf u}\in S^2$ with $D({\bf u},\OL\omega_{20})$ having exactly four distinct dot products are ${\bf u}\in\OL\omega_{12}$ with $D({\bf u},\OL\omega_{20})=\{\frac {\pm \varphi\pm 1}{\sqrt {3\varphi+6}}\}$.
\end {lemma}
\begin {proof}
Let ${\bf u}=(u_1,u_2,u_3)\in S^2$ be any vector, $a:=\max\{\left|u_i\right| : i=\OL {1,3}\}$ (we have $a>0$), and ${\bf x}=\frac {1}{a}\(\left|u_1\right|,\left|u_2\right|,\left|u_3\right|\)$. Then the set $D({\bf u},\OL\omega_{20})$, where $\OL\omega_{20}=\frac {1}{\sqrt {3}} B$, equals the set of dot products vector ${\bf x}$ forms with vectors from $\frac {a}{\sqrt{3}}B$. By Lemma~\ref {dodecahedrn}, vector ${\bf u}$ must form at least four distinct dot products with vectors from $\OL \omega_{20}$. It forms exactly four distinct dot products if and only if ${\bf x}$ is in $A_1$, which happens if and only if ${\bf u}\in \OL\omega_{12}$. 
Whenever ${\bf u}\in \OL \omega_{12}$, we have $a=\frac{\varphi}{\sqrt{\varphi+2}}$, and by Lemma \ref {dodecahedrn}, $D({\bf u},\OL\omega_{20})=\frac{\varphi}{\sqrt {3\varphi+6}}\{\pm 1\pm 1/\varphi\}=\{\frac {\pm\varphi\pm 1}{\sqrt {3\varphi+6}}\}$.
\end {proof}

\begin {proof}[Proof of Theorem \ref {dodecahedron}]The dodecahedron $\OL\omega_{20}=\{{\bf x}_1,\ldots,{\bf x}_{20}\}$ is an antipodal configuration on $S^2$ whose index set $\mathcal I_2(\OL\omega_{20})$ contains the set $\{1,2,3,4,5,7,8\}$, see Table 1. By Lemma \ref {ddcd}, any point ${\bf x}^\ast\in \OL\omega_{12}$ has $m=4$ distinct dot products with vectors from $\OL\omega_{20}$, which are $t_1=-\frac {\varphi+1}{\sqrt {3\varphi+6}}$, $t_2=\frac {-\varphi+1}{\sqrt {3\varphi+6}}$, $t_3=\frac {\varphi-1}{\sqrt {3\varphi+6}}$, $t_4=\frac {\varphi+ 1}{\sqrt {3\varphi+6}}$. Observe that
$$
\sum\limits_{i=1}^{4}t_i^2=\frac {2(\varphi+1)^2}{3\varphi+6}+\frac {2(\varphi-1)^2}{3\varphi+6}=\frac {4}{3}<\frac {28}{15}=\frac {m(2m-1)}{4m+d-3},
$$
where $m=4$. In view of Remark {3.2}, conditions \eqref {sumsq} hold. Then by Theorem~\ref {generalth},
$$
\min\limits_{{\bf x}\in S^2}\sum\limits_{i=1}^{20}f({\bf x}\cdot {\bf x}_i)=\sum\limits_{i=1}^{20}f({\bf x}^\ast\cdot {\bf x}_i)
$$
for any point ${\bf x}^\ast\in \OL\omega_{12}$.

Lemma \ref {ddcd} implies that for every ${\bf x}\in S^2\setminus \OL\omega_{12}$, the set $D({\bf x},\OL\omega_{20})$ is not contained in $\{\frac {\pm \varphi\pm 1}{\sqrt {3\varphi+6}}\}$. Then by Theorem \ref {generalth}, if $f^{(8)}>0$ on $(-1,1)$, point ${\bf x}$ is not a point of absolute minimum of $p_f(\cdot;\OL\omega_{20})$ over $S^2$.
\end {proof}

\section {Proof of Theorem \ref {E_8main}}\label {240}

A significant ingredient of the proof of Theorem \ref {E_8main} is the following lemma.

\begin {lemma}\label {E_8}
For every vector ${\bf x}\in S^7$, the set $D({\bf x},\OL\omega_{240})$ has at least five distinct elements. The only vectors ${\bf x}\in S^7$ such that $D({\bf x},\OL\omega_{240})$ has exactly five distinct elements are those in $\OL\omega_{2160}$ (the dot products are $\pm \frac {1}{2\sqrt {2}}$, $\pm \frac {1}{\sqrt {2}}$, and $0$) and those in $\OL \omega_{240}$ itself (the dot products are $\pm 1$, $\pm 1/2$, and~$0$).
\end {lemma}

\begin {proof}
We denote by $\W\omega_{112}$ the set of $112$ vectors on $S^7$ with $6$ zero coordinates and two coordinates with $\pm 1/\sqrt {2}$, and by $\W\omega_{128}$ the set of $128$ vectors with all eight coordinates $\pm \frac {1}{2\sqrt {2}}$ and an even number of ``$-$" signs. We have $\W\omega_{112}\cup\W\omega_{128}=\OL\omega_{240}$.
Let ${\bf x}=(x_1,\ldots,x_8)$ be an arbitrary vector on $S^7$, which forms at most five distinct dot products with vectors from $\OL\omega_{240}$. 

Assume to the contrary that the set $T:=\{\left|x_i\right|: i=\OL {1,8}\}$ has at least three distinct elements, and denote by $0\leq c< b< a$ some three of them. Then $D({\bf x},\W\omega_{112})$ contains at least $6$ distinct dot products: $\pm\(\frac {a}{\sqrt {2}}+ \frac {b}{\sqrt {2}}\)$, $\pm \(\frac {a}{\sqrt {2}}+ \frac {c}{\sqrt {2}}\)$, and $\pm \(\frac {b}{\sqrt {2}}+ \frac {c}{\sqrt {2}}\)$ contradicting the original choice of ${\bf x}$. Therefore, $T$ contains at most two distinct numbers.

Assume first that $T$ contains exactly two distinct elements and denote them by $0\leq b<a$. Let $b>0$. If $\left|x_i\right|=a$ for more than one index $i$, then the following $6$ distinct dot products are in $D({\bf x},\W\omega_{112})$: $\pm \sqrt {2} a$, $\pm\frac {a}{\sqrt {2}}\pm \frac {b}{\sqrt {2}}$. Therefore, $\left|x_i\right|=a$ for only one index $i$. Then $D({\bf x},\W\omega_{112})$ contains dot products $\pm \sqrt {2} b$, $ \pm\frac {a}{\sqrt {2}}\pm \frac {b}{\sqrt {2}}$, and $0$. They are all distinct if $a\neq 3b$ contradicting the original choice of ${\bf x}$. Therefore, $a=3b$. Since $\left|{\bf x}\right|=1$, we have $a=3/4$ and $b=1/4$. Then $D({\bf x},\W \omega_{112})$ contains $5$ distinct dot products: $\pm \frac{1}{\sqrt {2}}$, $\pm \frac {1}{2\sqrt {2}}$, and $0$. If ${\bf x}$ had an even number of negative components, we would let ${\bf z}\in\W \omega_{128}$ be the vector whose components had the same sign as the corresponding components of ${\bf x}$. Then ${\bf x}\cdot {\bf z}=\frac {5}{4\sqrt {2}}$, which would be the sixth dot product in $D({\bf x},\OL\omega_{240})$ contradicting the original choice of~${\bf x}$. Therefore, ${\bf x}$ has an odd number of negative components; that is ${\bf x}\in \OL\omega_{2160}$.

Now, assume that $b=0$. Let $k$ be the number of non-zero coordinates of ${\bf x}$. Then $a=\frac {1}{\sqrt {k}}$. If $k=3,5, 6$, or $7$, it is not difficult to see that the set $D({\bf x},\W \omega_{112})$ contains dot products $0,\pm \frac {1}{\sqrt{2k}},\pm \frac {2}{\sqrt{2k}}$. With ${\bf z}=\(\frac {1}{2\sqrt{2}},\ldots,\frac {1}{2\sqrt{2}}\)\in \W\omega_{128}$ we would get the sixth dot product ${\bf x}\cdot {\bf z}=\frac {k}{2\sqrt{2k}}$. Therefore, $k$ is $1$, $2$, or $4$. If $k=4$, then four coordinates of ${\bf x}$ equal $\pm 1/2$ and the remaining four equal~$0$; that is, ${\bf x}\in \OL\omega_{2160}$. If $k=1$, then ${\bf x}$ has seven zero coordinates with the remaining coordinate equal $\pm 1$ which puts ${\bf x}$ in $\OL\omega_{2160}$ again. When $k=2$, we have ${\bf x}\in \OL\omega_{240}$.

Finally, assume that $T$ contains only one element (denote it by $a$). Then every coordinate of ${\bf x}$ is $\pm \frac {1}{2\sqrt {2}}$ and $D({\bf x},\W\omega_{112})$ contains dot products $\pm 1/2,0$. Let $i$ be the number of negative coordinates in ${\bf x}$. Assume to the contrary that $i$ is odd. We take ${\bf z}\in \W\omega_{128}$ with first seven coordinates having the same sign as the corresponding coordinates of ${\bf x}$. Then the eighth coordinates of ${\bf z}$ and ${\bf x}$ have opposite signs; that is, ${\bf x}\cdot {\bf z}=3/4$. We next change the sign of first two, four, or six coordinates of ${\bf z}$. Then ${\bf z}$ still has an even number of negative coordinates with ${\bf x}\cdot {\bf z}=1/4$, $-1/4$, or $-3/4$. Thus, when $i$ is odd, $D({\bf x},\OL\omega_{240})\supset \{0,\pm 1/4,\pm 1/2,\pm 3/4\}$ contradicting the original choice of ${\bf x}$. Therefore, $i$ is even and ${\bf x}\in \OL\omega_{240}$.

Thus, if ${\bf x}\notin \OL\omega_{240}\cup \OL\omega_{2160}$, the set $D({\bf x},\OL\omega_{240})$ contains more than five distinct dot products. If ${\bf x}\in \OL \omega_{240}$ then $D({\bf x},\OL\omega_{240})=\{0,\pm 1/2,\pm 1\}$, and if ${\bf x}\in \OL\omega_{2160}$ then $D({\bf x},\OL\omega_{240})=\{\pm \frac{1}{2\sqrt{2}},\pm \frac{1}{\sqrt{2}},0\}$.
  \end {proof}

\begin {proof}[Proof of Theorem \ref {E_8main}]
The $E_8$ lattice $\OL\omega_{240}=\{{\bf x}_1,\ldots,{\bf x}_{240}\}$ is an antipodal configuration on $S^7$ whose index set $\mathcal I_7(\OL\omega_{240})$ contains the set $\{1,2,3,4,5,6,7,9,10\}$, see Table 1. By Lemma \ref {E_8}, any point ${\bf x}^\ast\in \OL\omega_{2160}$ has $m=5$ distinct dot products with vectors from $\OL\omega_{240}$, which are $t_1=-\frac {1}{\sqrt {2}}$, $t_2=-\frac {1}{2\sqrt {2}}$, $t_3=0$, $t_4=\frac {1}{2\sqrt {2}}$, $t_5=\frac {1}{\sqrt {2}}$. Observe that
$$
\sum\limits_{i=1}^{5}t_i^2=\frac {5}{4}<\frac {15}{8}=\frac {m(2m-1)}{4m+d-3},
$$
where $m=5$. In view of Remark \ref {3.2}, conditions \eqref {sumsq} hold. Then by Theorem~\ref {generalth},
$$
\min\limits_{{\bf x}\in S^7}\sum\limits_{i=1}^{240}f({\bf x}\cdot {\bf x}_i)=\sum\limits_{i=1}^{240}f({\bf x}^\ast\cdot {\bf x}_i)
$$
for any point ${\bf x}^\ast\in \OL\omega_{2160}$.

By Lemma \ref {E_8}, for every ${\bf x}\in S^7\setminus \OL\omega_{2160}$, the set $D({\bf x},\OL\omega_{240})$ is not contained in $\{\pm \frac {1}{2\sqrt {2}},\pm \frac {1}{\sqrt {2}},0\}$. Then by Theorem \ref {generalth}, if $f^{(10)}>0$ on $(-1,1)$, point ${\bf x}$ is not a point of absolute minimum of $p_f(\cdot;\OL\omega_{240})$ over $S^7$.
\end {proof}

}

\begin {thebibliography}{99}
\bibitem {AbrSte1965}
M.Abramowitz, I.A.Stegun, {\it Handbook of Mathematical Functions}, Dover, New York, 1965.
\bibitem {Bil2015}
M. Bilogliadov, Equilibria of Riesz potentials generated by point charges at the roots of unity, {\it Comput. Methods Funct. Theory} {\bf 15} (2015), no. 4, 471--491.
\bibitem {Borstiff}
S.V. Borodachov, Absolute minima of potentials of a certain class of spherical designs (in preparation).
\bibitem{Bor2022talk}
S.V. Borodachov, Min-max polarization for certain classes of sharp configurations on the sphere, {\it Workshop "Optimal Point Configurations on Manifolds"}, ESI, Vienna, January 17--21, 2022. https://www.youtube.com/watch?v=L-szPTFMsX8
\bibitem {BortalkGMU}
S.V. Borodachov, Extreme values of potentials of spherical designs and the polarization problem, {\it XII International Conference of the Georgian Mathematical Union}, August 29--September 3, 2022, Batumi, Georgia.
\bibitem {Borsimplex}
S.V. Borodachov, Polarization problem on a higher-dimensional sphere for a simplex, {\it Discrete and Computational Geometry} {\bf 67} (2022), no. 2, 525--542.
\bibitem {BorMinMax}
S.V. Borodachov, Min-max polarization for certain classes of sharp configurations on the sphere (submitted). https://arxiv.org/pdf/2203.13756.pdf.
\bibitem{BorHarSafbook}
S. Borodachov, D. Hardin, E. Saff, {\it Discrete Energy on Rectifiable Sets}. Springer, 2019.
\bibitem {BoyDanKaz2001}
P. Boyvalenkov, D. Danev, P. Kazakov, Indexes of spherical codes, {\it Codes and association schemes} (Piscataway, NJ, 1999), 47–57. {\it DIMACS Ser. Discrete Math. Theoret. Comput. Sci.}, {\bf 56}, Amer. Math. Soc., Providence, RI, 2001.
\bibitem {BoyDraHarSafSto600cell}
P.G. Boyvalenkov, P.D. Dragnev, D.P. Hardin, E.B. Saff, M.M. Stoyanova, On polarization of spherical codes and designs (submitted). https://arxiv.org/abs/2207.08807.
\bibitem{CohKum2007}
H. Cohn, A. Kumar, Universally optimal distribution of points on spheres, {\it J. Amer. Math. Soc.} {\bf 20} (2007), no. 1, 99--148.
\bibitem {ConSlo1999}
J. Conway, N.J.A. Sloane, {\it Sphere packings, lattices, and groups}, Springe Science $+$ Business Media, 1999. New York. 
\bibitem {DelGoeSei1977}
P. Delsarte, J.M. Goethals, J.J. Seidel, Spherical codes and designs, {\it Geometriae Dedicata}, {\bf 6} (1977), no. 3, 363--388.
\bibitem {GioKhi2020}
G. Giorgadze, G. Khimshiashvili, Stable equilibria of three constrained unit charges, {\it Proc. I. Vekua Inst. Appl. Math.} {\bf 70} (2020), 25--31.
\bibitem {IsaKel1965}
E. Isaacson, H.B. Keller, {\it Analysis of Numerical Methods}, Dover Books on Mathematics, 1994.
\bibitem {Lev1979}
V.I.Levenshtein, On bounds for packings in $n$-dimensional Euclidean space,
{\it Soviet Math. Dokladi}, 20, 1979, 417--421.
\bibitem {Lev1992}
V.I. Levenshtein, Designs as maximum codes in polynomial metric spaces,
{\it Acta Appl. Math.} {\bf 25} (1992), 1--82.
\bibitem {Lev1998}
V.I. Levenshtein, Universal bounds for codes and designs, Chapter 6 in
{\it Handbook of Coding Theory}, V. Pless and W.C. Huffman, Eds., Elsevier Science B.V., 1998.
\bibitem {Mus2003}
O.R. Musin, The kissing number in four dimensions, {\it Annals of Mathematics}, {\bf 168} (2008), 1--32.
\bibitem {NikRaf2011}
N. Nikolov, R. Rafailov, On the sum of powered distances to certain sets of points on the circle, {\it Pacific J. Math.} {\bf 253} (2011), no. 1, 157--168.
\bibitem{NikRaf2013}
N. Nikolov, R. Rafailov, On extremums of sums of powered distances to a finite set of points. {\it Geom. Dedicata} {\bf 167} (2013), 69--89.
\bibitem {Sto1975circle}
K. Stolarsky, The sum of the distances to certain pointsets on the unit circle, {\it Pacific J. Math.} {\bf 59} (1975), no. 1, 241--251. 
\bibitem {Sto1975}
K. Stolarsky, The sum of the distances to $N$ points on a sphere, {\it Pacific J. Math.} {\bf 57} (1975), no. 2, 563--573.
\end {thebibliography}
\end {document}